\documentclass[oneside,a4paper,11pt,reqno]{amsart}
\usepackage{array,amsfonts,amsmath,amsthm}
\usepackage[dvips]{graphicx}
\usepackage{a4wide}
\usepackage[usenames]{color}
\definecolor{webgreen}{rgb}{0,.5,0}
\definecolor{webbrown}{rgb}{.6,0,0}
\usepackage[colorlinks=true,
	linkcolor=webgreen,
	filecolor=webbrown,
	citecolor=webgreen]{hyperref}

\newtheorem{theorem}{Theorem}

\newtheorem{corollary}[theorem]{Corollary}

\newtheorem{fact}[theorem]{Fact}
\theoremstyle{definition}

\theoremstyle{remark}
\newtheorem{remark}{Remark}
\newtheorem{example}{Example}

\newcommand{\weight}{\mathcal{W}}
\newcommand{\property}{\emph{property}}
\newcommand{\dyck}{\mathbb{D}}

\newcommand{\seqnum}[1]{\href{http://www.research.att.com/~njas/sequences/#1}{#1}}

\begin{document}

\title[Parametric production matrices and weighted succession rules]{Parametric production matrices and weighted succession rules: a Dyck path example}
\author[R. Parviainen]{Robert Parviainen}
\address{ARC Centre of Excellence for Mathematics 
   	and Statistics of Complex Systems\\
  	139 Barry Street, The University of Melbourne, Victoria, 3010}
\email{robertp@ms.unimelb.edu.au}

\date{\today}
\begin{abstract}
We introduce weighted succession rules and parametric production matrices ---  simple extensions of the standard ECO method succession rules and production matrices. The purpose is to enumerate combinatorial objects with respect to several variables. We consider one main example, from the theory of Dyck paths. The path statistics primarily considered are peak height, rise height, and certain subwalks. Many classical sequences (such as the Catalan, Motzkin, Narayana and Schroeder numbers) are incorporated in this example. 
\end{abstract}
\maketitle

%%%
%%%
%%%

\section{Introduction}
Recently Deutsch, Ferrari and Rinaldi, \cite{DFR2005}, introduced production matrices as a device for ECO method calculations. Building on the concept of generating trees, the \emph{Enumeration of Combinatorial Objects} method was developed and systematised by Barcucci and collaborators, \cite{BDLPP1999}. It is a general method for enumerating combinatorial objects by studying how the objects grow according to a parameter. At the centre of method are the succession rules that describe how the objects grow --- production matrices translate these to matrix notation.

The production matrices considered by Deutsch et al. are integer valued. The step to include matrices with parameters is natural and straightforward, and corresponds to weighted succession rules. In this paper we study a rich example of a parametric production matrix. Although only one matrix, it embodies many well known combinatorial number sequences, such as the Catalan, Motzkin, and Narayana numbers. The underlying objects for our matrix are Dyck paths, and the main statistics considered are heights of peaks and rises, and occurrences of certain subwalks. We recover known results about the height distribution of peaks, and use the matrix to show that it is equivalent to a rise height distribution. 

In Section \ref{section:definitions} we give a brief description of the ECO method, define succession rules and production matrices, and introduce weighted successions.  We recall results from \cite{DFR2005} for operations on product matrices (which continue to hold in the parametric case).  In Section \ref{section:main} we define Dyck paths and some of their statistics. We then define a parametric production matrix and show how it is related to those statistics. 

%%%
%%%
%%%

\section{The ECO method}\label{section:definitions}
The core of the ECO method is a recursive description of a class of combinatorial objects. This should be done in such a way that, if $\mathcal O_{n}$ denotes the set of objects of size $n$, each object $O^{\prime}\in\mathcal O_{n+1}$ is achieved from one and only one object $O\in\mathcal O_{n}$. We say that $O^{\prime}$ is a \emph{successor} of $O$. 

We assign a label $(k)$, $k\in\mathbb{N^{+}}$ (the positive integers), to each object. An object's label give the number of successors of that object. The succession rule dictates the labels of these successors. The rule also includes an axiom $(a)$, $a\in\mathbb{N^{+}}$, which gives the label of the smallest object. In the basic case a succession rule $\Omega$ is written as
\[
	\Omega=
	\begin{cases}
		(a)\\
		(k)\leadsto\left(e_{1}(k)\right)\left(e_{2}(k)\right)\cdots\left(e_{k}(k)\right)
	\end{cases},
\]
where $e_{i}(k):\mathbb{N^{+}}\to\mathbb{N^{+}}$ gives the labels of the $k$ successors of an object with label $k$.

It is natural to think of a succession rule as a generating tree: the root of the tree is the axiom, $(a)$, and if a node has label $(k)$, it has $k$ children labelled $e_{1}(k), e_{2}(k), \ldots, e_{k}(k)$.

A succession rule $\Omega$ defines a sequence of positive integers $a_{n}=$ the number of nodes at level $n$, $n\geq 0$. The generating function for $a_{n}$ is defined by $f_{\Omega}(z)=\sum_{n}a_{n}z^{n}$.

The production matrix $P(\Omega)=P=(p_{i,j})_{i,j\geq 1}$ of a succession rule is defined as follows. List all labels of a rule. Let this list be $l_{1}, l_{2}, \ldots$. Then $p_{i,j}$ equals the number of successors with label $l_{j}$ produced by an object with label $l_{i}$.

\begin{example}
The succession rule 
\[
	\Omega=
	\begin{cases}
		(1)\\
		(1)\leadsto(2)\\
		(2)\leadsto(1)(2)\\
	\end{cases},
\]
has production matrix
\[
	P(\Omega)=\left(
	\begin{array}{cc}
		0 & 1 \\
		1 & 1
	\end{array}\right),
\]
since label $(1)$ has one successor with label $(2)$ and $(1)$ has two successors with labels $(1)$ and $(2)$. The sequence defined by $\Omega$ is the Fibonacci sequence (A000045\footnote{Six-digit numbers
prefixed by `A' indicate the corresponding entry in \emph{The On-Line Encyclopedia of Integer Sequences}
 \cite{OEIS}.}), $a_{n}=1,1,2,3,5,8, \ldots$.
\end{example}

\subsection{Parametric production matrices and weighted succession rules}
Naturally, one would like to use the ECO method to obtain not only the size distribution of the objects, but also multi--dimensional distributions. 

One strategy is to choose the labels wisely, since the distribution of labels on objects of a given size might be acquired. Indeed, using production matrices this is essentially as easy (or hard) as getting the size distribution. Another strategy is to use the succession rule and other data to find functional relations for the multivariate generating function.

One of the aims of this paper is to demonstrate an alternative way of using production matrices to find multivariate distributions, by introducing weighted succession rules and parametric production matrices. The idea is simple: 

Assume we are interested in the distribution of \property, jointly with size. If a succession $(k)\leadsto (l)$ adds one \property, we give that succession weight $x$. The sequence generated by the succession rule now is the weighted sum over objects of given size, and the generating function for the sequence is the bivariate generating function for \property~and size.

Let $\weight$ be a set of weights, and consider a matrix $P=(p_{i,j})_{i,j\geq 1}$ where $p_{i,j}\in\weight$. This \emph{parametric} production matrix corresponds to a \emph{weighted} succession rule. A weighted succession rule is written as 
\[
	\Omega=
	\begin{cases}
		\left(a; w\left(a\right)\right)\\
		(k)\leadsto\left(e_{1}(k); w\left(e_{1}(k)\right)\right)\left(e_{2}(k); 
			w\left(e_{1}(k)\right)\right)\cdots\left(e_{k}(k); w\left(e_{1}(k)\right)\right)
	\end{cases},
\]
where again $e_{i}(k):\mathbb{N^{+}}\to\mathbb{N^{+}}$ gives the labels of the $k$ successors of an object with label $k$, and 
$w\left(e_{i}(k)\right)\in\weight$ gives the weight of that particular succession. %We also give the axiom $a$ a weight $w(a)$.
In terms of generating trees, each edge is given a weight from $\weight$. We also have the possibility to give the axiom $a$ a weight $w(a)$. Of course, if $\weight=\{1\}$ we recover the unweighted succession rules.

From a weighted succession rule $\Omega$, a parametric production matrix $P(\Omega)=P=(p_{i,j})_{i,j\geq 1}$ is defined as follows. List all labels of a rule. Let this list be $l_{1}, l_{2}, \ldots$. Then $p_{i,j}$ equals the sum of weights of successors with label $l_{j}$ produced by an object with label $l_{i}$. The sequence $\{a_{n}\}_{n\geq 0}$ produced by $\Omega$ is now a sequence of polynomials in $\weight$, and we write
\[
	f_{\Omega}(\weight, z)=\sum_{n=0}^{\infty}a_{n}z^{n}.
\]
%=\sum_{n=0}^{\infty}\sum_{k_{1}, k_{2}, \ldots} w_{1}^{k_{1}}w_{2}^{k_{2}}\cdots z^{n}.
%\]

\begin{example}
The succession rule 
\[
	\Omega=
	\begin{cases}
		(1;1)\\
		(1)\leadsto(1;1)\\
		(3)\leadsto(1;x)(2;1)\\
	\end{cases},
\]
has production matrix
\[
	P(\Omega)=\left(\begin{array}{cc}0 & 1 \\ x & 1\end{array}\right),
\]
The sequence defined by $\Omega$ gives the Fibonacci polynomials (\seqnum{A011973}) $a_{n}=1, 1, 1+x, 1+2x, 1+3x+x^{2}, 1+4x+3x^{2}, \ldots$.
\end{example}

Let us fix some notation. Let $u^{\top}=\{1\ 0\ 0\ \ldots\}$,  $e=\{1\ 1\ 1\ \ldots \}^{\top}$ and $I$ define a row vector, column vector, and identity matrix, respectively, of appropriate size. If a succession rule $\Omega$ defines a sequence $a_{n}$ with generating function $f_{\Omega}(\weight, z)$, we say that $P=P(\Omega)$ defines $a_{n}$ and write $P\to a_{n}$ and $P\to f_{\Omega}(\weight, z)$.

We now cite some useful results for production matrices. These are from \cite{DFR2005}, where they were stated in terms of integer valued production matrices. To verify that these still hold in the parametric case is trivial.
\begin{fact}[Fact (iii), \cite{DFR2005}]\label{fact:seq}
	If $P\to a_{n}$ then $a_{n}=u^{\top}P^{n}e$.
\end{fact}
\begin{fact}[Fact (v), \cite{DFR2005}]\label{fact:gf}
	If $P\to f(\weight, z)$ then $f(\weight, z)=u^{\top}(I-z P)^{-1}e$.
\end{fact}

\begin{theorem}[Theorem 3.2, \cite{DFR2005}]\label{theorem:gf}
If $b$, $c$ and $r$ are nonnegative integers or weights and $P\to f(\weight, z)$, then
\[
	M=
	\left(
	\begin{array}{cc}
		b & ru^\top \\
		ce & P
	\end{array}
	\right)\to\frac{1+rzf(\weight, z)}{1-bz-rcz^{2}f(\weight, z)}.
\]
\end{theorem}
\begin{corollary}[Corollary 3.1, \cite{DFR2005}]\label{corollary:recursion}
If $b$, $c$ and $r$ are nonnegative integers or weights and 
\[
	P=
	\left(
	\begin{array}{cc}
		b & ru^\top \\
		ce & P
	\end{array}
	\right)\to f(\weight, z),
\]
then $f(\weight, z)$ satisfies 
\[
	rcz^{2}f^{2}(\weight,z)-(1-(b+r)z)f(\weight, z)+1=0.
\]
\end{corollary}

%%%
%%%
%%%

\section{The Main Example}\label{section:main}
\subsection{Dyck paths}
A Dyck path is a lattice path in the first quadrant from $(0,0)$ to $(2n,0)$ that uses steps $(1,1)$, we call these up, or $u$, steps, and $(1,-1)$ (down, or $d$, steps). We will call $n$ the \emph{semilength} of the path. Let $\dyck_{n}$ denote the set of Dyck paths of semilength $n$. The number of Dyck paths of semilength $n$ are given by the $n$th Catalan number $c_{n}$ (\seqnum{A000108}), 
\[
	c_{n}=\frac{1}{n+1}\binom{2n}{n}.
\]
The generating function for the Catalan numbers satisfies
\[
	C(z)=1+z C^{2}(z)=\sum_{n=0}^{\infty}c_{n}z^{n}=\frac{1-\sqrt{1-4z}}{2z}.
\]
%so the generating function for Dyck paths is $C^{2}(z)$.

A \emph{peak} in a Dyck path $p$ is an occurrence of an up step immediately followed by a down step. The {height} of the peak is the height of the start point of the $u$ step.
A \emph{high peak} is a peak at height 1 or higher.
A \emph{rise}  is a maximal sequence of consecutive $u$ steps. The {height} of the rise is the number of $u$ steps minus 1.
A \emph{high rise} is a rise with positive height.
A \emph{contact} is any place the path touches the $x$-axis. 
A \emph{segment} is any subpath consisting of consecutive steps. The segment is a \emph{level} segment if it starts and ends at the same height.
An \emph{excursion} is the segment of the walk between two contacts with no other contact in between. 
Note that a Dyck path $p$ with $k+1$ contacts can be decomposed uniquely into its $k$ excursions $p=e_{1}e_{2}\cdots e_{k}$.

\begin{figure}[htbp]
\begin{center}
\includegraphics[scale=1.2]{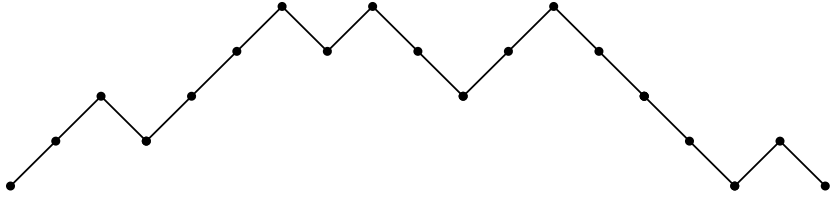}
\caption{A Dyck path of semilength 9. The height of the rises are, in order, 1,2,0,1,0. The path has 4 high peaks, of heights 1,3,3,3. Using  succession rule $\Omega$ of Theorem \ref{theorem:main} the path will have weight $x_{0}^{5}x_{1}x_{2}^{3}y_{1}^{1}y_{2}^{2}y_{4}^{1}$.  Using succession rule $\Omega^{\prime}$ of Theorem \ref{theorem:equi} the path will have weight $y_{1}^{1}y_{3}^{3}$.}
\label{fig:path}
\end{center}
\end{figure}

\subsection{The matrix}
We are interested in the following matrix.
\begin{equation*}
	P=
	\left(
	\begin{array}{ccccccc}
		x_{1}y_1 & x_{0} & 0 & 0 & 0 & 0 & \cdots \\
		x_{2}y_1 & x_{1}y_2 & x_{0} & 0 & 0 & 0 &  \cdots \\
		x_{3}y_1 & x_{2}y_2 & x_{1}y_3 & x_{0} & 0 & 0 &  \cdots \\
		x_{4}y_1 & x_{3}y_2 & x_{2}y_3 & x_{1}y_4 & x_{0} & 0 &  \cdots \\
		x_{5}y_1 & x_{4}y_2 & x_{3}y_3 & x_{2}y_4 & x_{1}y_5 & x_{0} &  \cdots \\
		x_{6}y_1 & x_{5}y_2 & x_{4}y_3 & x_{3}y_4 & x_{2}y_5 & x_{1}y_6 &  \cdots \\
		\vdots & \vdots & \vdots & \vdots & \vdots & \vdots & \ddots
	\end{array}
	\right).
\end{equation*}
It is, as is shown below, the production matrix for Dyck paths, where the $y_{k}$ marks the height of rise number $k$ and $x_{k}$ marks level segments of the form $u(uq_{k}d)^{k}d$, where $q_{k}$ are Dyck paths. We note that many classical number sequences arise from specialisations of the sequence associated with $P$. See Section \ref{section:spec}.

\begin{theorem}\label{theorem:main}
For a Dyck path $p$, let $h_{k}(p)$ denote the height of rise number $k$ in $p$ and let $s_{k}(p)$ denote the number of level segments of type $u(uq_{k}d)^{k}d$, where $q_{k}$ are Dyck paths. Then the weighted sum $a_{n}$ over Dyck paths of semilength $n$ is given by
\[
	a_{n}
	:=\sum_{p\in\dyck_{n}}x_{0}^{h_{0}(p)}y_{1}^{s_{1}(p)}x_{1}^{h_{1}(p)}y_{2}^{s_{2}(p)}\cdots
	=u^{\top}P^{n}e,
\]
and the generating function is given by
\begin{align*}
	F(x_{0}, x_{1},\ldots; y_{1}, y_{2}, \ldots; z)
	:&=\sum_{n}z^{n}\sum_{p\in\dyck_{n}}x_{0}^{h_{0}(p)}y_{1}^{s_{1}(p)}x_{1}^{h_{1}(p)}y_{2}^{s_{2}(p)}\cdots \\
	&=u^{\top}(I-z P)^{-1}e.
\end{align*}
\end{theorem}
\begin{figure}[htbp]
\begin{center}
\includegraphics[scale=0.5, angle=270]{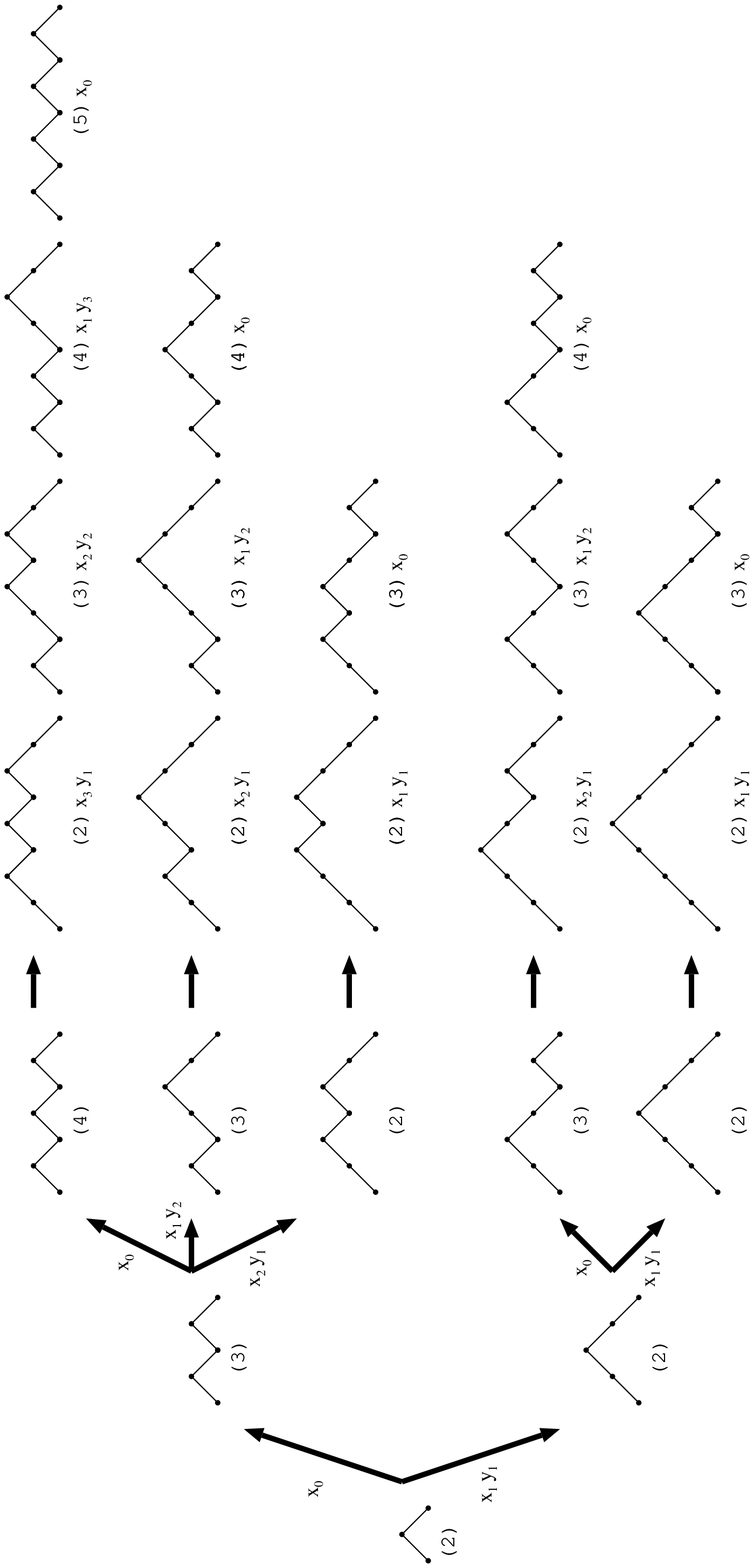}
\caption{The generating tree associated with the succession rule of Theorem \ref{theorem:main}.}
\label{fig:gen tree 1}
\end{center}
\end{figure}

\begin{proof}
We will describe a weighted succession rule $\Omega$ for Dyck paths whose parametric production matrix is $P$. The labels are the number of contacts (= the number of excursions + 1).  The axiom is the path $ud$ with label $(2; x_{0})$.
The generating tree for this rule is sketched in Figure \ref{fig:gen tree 1}.

Given a Dyck path $p$ of length $n$ with $k$ excursions, $p=e_{1}e_{2}\cdots e_{k}$ we now explain how the successors are derived. There are $k+1$ successors to $p$, with $1, 2, \cdots, k+1$ excursions. They are, respectively,
\begin{description}
	\item[1] $ue_{1}e_{2}\cdots e_{k}d$
	\item[2] $e_{1}ue_{2}\cdots e_{k}d$\\$\vdots$
	\item[k] $e_{1}e_{2}\cdots e_{k-1}ue_{k}d$
	\item[k+1] $e_{1}e_{2}\cdots e_{k}ud$.
\end{description}
Now, let a succession from a path with $k$ excursions to one with $l$, $1\leq l \leq k+1$, have weight $y_{l}x_{k-l+1}$. The succession increases the height of rise number $l$ by one and adds one excursion, by default a level segment, of type $u(uqd)^{k-l+1}d$ to the path as desired. Furthermore, each path has exactly one parent: the path obtained by by lowering the last excursion $uqd$ to $q$, a Dyck path. 
\end{proof}

In the following, we will use the notation $x=(x_{0}, x_{1}, x_{2}, \ldots)$ and  $y=(y_{1}, y_{2}, y_{3}, \ldots)$. If $x_{k}=1$ for $k>n$ we will write $F(x_{1}, \ldots, x_{n}; y;z)$ instead of $F(x; y; z)$, and similarly for $y$.
\begin{corollary}\label{cor:f0}
\[
	F_{0}:=F(x_{0}, x_{1};; z) = \frac{1 - (x_0 + x_1) z -\sqrt{ (1 - (x_0 +x_1) z)^2-4 x_0 z^2}}{2 x_0 z^2}
\]
\end{corollary}
\begin{proof}
It follows from Corollary \ref{corollary:recursion} that 
\[
	x_{0}z^{2}F_{0}^{2}-\left(1-(x_{0}+x_{1})z\right)F_{0}+1=0.
\]
\end{proof}
Applying Theorem \ref{theorem:gf} we get
\begin{corollary}\label{cor:f1}
\[
	F_{1}:=F(x_{0}, x_{1}; y_{1}; z) = 
	\frac{1+x_{0} z F_{0}(x_{0},x_{1};;z)}{1-x_{1}y_{1}z-x_{0}y_{1}z^{2}F_{0}(x_{0},x_{1};;z)},
\]
and in general
\[
	F_{k}:=F(x_{0}, x_{1}; y_{1}, y_{2}, \ldots, y_{k}; z) = 
	\frac{1+x_{0} z F_{k-1}(x_{0},x_{1}; y_{2},\ldots, y_{k};z)}
	{1-x_{1}y_{1}z-x_{0}y_{1}z^{2}F_{k-1}(x_{0},x_{1}; y_{2},\ldots, y_{k};z)}.
\]
Note the use of $F_{k-1}(x_{0},x_{1}; y_{2},\ldots, y_{k};z)$ and not of $F_{1}(x_{0},x_{1}; y_{1},\ldots, y_{k-1};z)$.
\end{corollary}

%%
%%

%\subsection{Excursions}

%%
%%

%\subsection{Rise and peak heights}
Below we find bivariate generating functions $F(1,\ldots, 1,x_{k};1;z)$ for length and height of rise number $k$ in a Dyck path. First we show that the heights of high rises have the same distribution as the number of high peaks. The reason for not using the perhaps more familiar high peak statistics in the definition of the main production matrix is that the \emph{joint} distribution of number of high peaks and level segments is not the same as the joint distribution of heights of rises and level segments. 

\begin{figure}[htbp]
\begin{center}
\includegraphics[scale=0.75, angle=270]{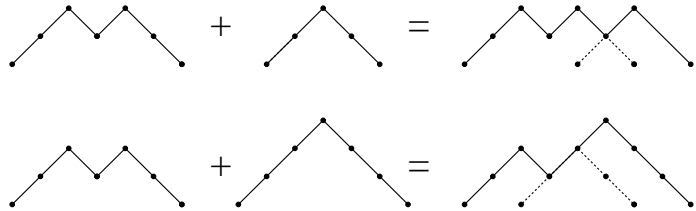}
\caption{Example of the use of the concatenation rule used in the proof of Theorem \ref{theorem:equi}.}
\label{fig:concat}
\end{center}
\end{figure}

\begin{theorem}\label{theorem:equi}
The number of Dyck paths with whose high rises have heights (in order) $k_{1}, \ldots, k_{l}$ equals the number of Dyck paths with $k_{1}, \ldots, k_{l}$ high peaks of heights $1, \ldots, l$, respectively.
\end{theorem}
\begin{figure}[htbp]
\begin{center}
\includegraphics[scale=0.5, angle=270]{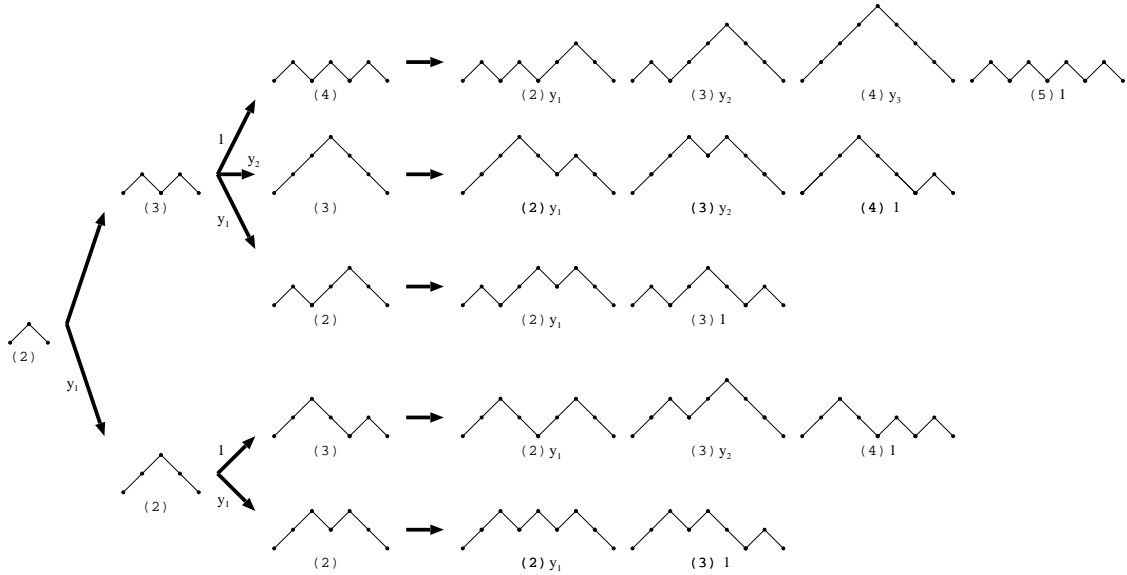}
\caption{The generating tree associated with the succession rule of Theorem \ref{theorem:equi}.}
\label{fig:gen tree 2}
\end{center}
\end{figure}
\begin{proof}
The result can be proved by describing an appropriate succession rule $\Omega^{\prime}$. This will be such that 
$P(\Omega^{\prime})=P(\Omega)|_{x=0}$.
The generating tree for this rule is sketched in Figure \ref{fig:gen tree 2}. 

Paths are labeled by the number of steps in the last descent minus 1, except if that number is 0 - then the label is the number of successive $ud$s at the end, which equals the parent's label plus 1. The axiom is the path $ud$, and is given weight 1.
The successors of a path $p$ with label $(k)$ are paths $p_{1}, \ldots, p_{k+1}$ with labels $(1)$ to $(k+1)$.
The path $p_{k+1}$ is $pud$. The path $p_{l}$, $l\leq k$, is achieved from $p$ as follows. It is a combination of the paths $p$ and $q=u^{l+1}d^{l+1}$;  the latter is placed such that it ends at $2n+2$ and in $p_{l}$ we follow $p$ up to the first intersection of $p$ and $q$, after which we follow $q$. See Figure \ref{fig:concat} for illustration.

Now, since a path's label gives either the number of $ud$s at the end of $p$, or 1 plus the number of $d$ steps in the last descent, the high peaks of $p$ are preserved in $p_{l}$, and if, $l\leq k$, $p_{l}$ has an extra peak at height $l$.

It is thus clear that if we give the updates of $(k)$ weights  $y_{1}, \ldots, y_{k}, 1$ the succession rule corresponds to heights of high peaks as desired. It is also easy to see that each path has a unique parent. 
\end{proof}

The next result is essentially equivalent to Theorem 1.1(ii), \cite{mans2002}. Mansour gives generating functions for $k$ peaks of height $h$ in terms of Chebyshev polynomials. These can be summed to give an equivalent function for peaks of height $h$.  As a function of $t$, the function $G_{n}$ below is of the form $\frac{a+t b}{c+d t}$, and so can easily be expanded in $t$ to give one--dimensional generating functions if desired.

\begin{corollary}
Let $G_{n}(t,z)$ be the generating function for high peaks of height $n\geq 1$,
\[
G_{n}(t,z)=\sum_{i=0}^{\infty}\sum_{j=0}^{\infty}g_{i,j}z^{i}t^{j},
\]
where $g_{i,j}$ is the number of Dyck paths of semilength $i+1$ with $j$ high peaks at height $n$. 

Then
\[
	G_{1}(t,z)=\frac{1+z C(z)^{2}}{1-t z-t z^{2}C(z)^{2}},
\]
and for $n\geq2$
\[
	G_{n}(t,z)=\frac{G_{1}(t,z) z U(n-4,\frac{1}{2\sqrt{z}})-U(n-2,\frac{1}{2\sqrt{z}})}
	{G_{1}(t,z) z^{2}U(n-2,\frac{1}{2\sqrt{z}})-zU(n,\frac{1}{2\sqrt{z}})}
\]
\end{corollary}

\begin{proof}
The first statement follows from Corollaries \ref{cor:f0} and \ref{cor:f1}.

Let $P_{n}$ be the matrix $P$ with all parameters except $y_{n}$ set to 1. Now
\[
	P_{n}=\left(
	\begin{array}{cc}
		1&u^{\top}\\
		e&P_{n-1}
	\end{array}
	\right).
\]
By Theorem \ref{theorem:gf} the generating function for $P_{n}, n\geq 2$ is 
\[
	G_{n}=\frac{1+z G_{n-1}}{1-z-z^{2}G_{n-1}}
\]
Now, straightforward but cumbersome calculations show that 
$\frac{G_{1} z U(n-4,\frac{1}{2\sqrt{z}})-U(n-2,\frac{1}{2\sqrt{z}})}
{G_{1} z^{2}U(n-2,\frac{1}{2\sqrt{z}})-U(n,\frac{1}{2\sqrt{z}})}$ satisfies the recursion. 
\end{proof}
\begin{remark}
The stated formula for $G_{n}$ in fact holds for $n=0$ and $n=1$ as well.
\end{remark}
\subsection{Examples}\label{section:spec}
As mentioned earlier, many well--known combinatorial sequences are contained in the sequence produced by the main matrix $P$. Here we list a selection of the most interesting of these. 
\begin{itemize}
\item $F( ; ;z)\to {1,2,5,14,\ldots}$ =
the Catalan numbers = the number of Dyck paths. (\seqnum{A000108}).
\item $F(1,0;;z)\to {1,1,2,4,9,\ldots}$ =
the Motzkin numbers = the number of Motzkin paths. (\seqnum{A001006}).
\item $F(2; ;z)\to 1,3,11,45,\ldots$ = the little Schroeder numbers =  the number of Schroeder paths with no peaks at height 1 (\seqnum{A001003}).
\item $F(s;;z)\to {1, 1+t, 1+3t+t^{2}, 1+6t+6t^{2}+t^{3},\ldots}$ =
the Narayana numbers = the number of Dyck paths with $k$ peaks. (\seqnum{A001263}).
\item $F(2s;;z)\to {1, 1+2t, 1+6t+4t^{2}, 1+12t+24t^{2}+8t^{3},\ldots}$ =
the number of double rise bi-colored  Dyck paths with $k$ double rises. (\seqnum{A114687}).
\item $F(1+s;;z)\to {1, 2+t, 5+5t+t^{2}, 14+21t+9t^{2}+t^{3},\ldots}$ =
the number of Schroeder paths with $k$ peaks, but none at height 1 =
the number of double rise bi-colored  Dyck paths with $k$ double rises of a given color. (\seqnum{A126216}).
\item $F(1,2;;z)\to {1,3,10,36,137,\ldots}$ =
the number of Schroeder paths with 0 peaks at odd level. (\seqnum{A002212}).
\item $F(1,s;;z)\to {1, 1+t, 2+2t+t^{3}, 4+6t+3t^{2}+t^{3}, \ldots}$ =
the number of Dyck paths with $k$ peaks at even height. (\seqnum{A091869}).
\item $F(t, 2; ;z)\to  1, 2+t,  4+ 5t+t^{2}, 8+18t+9t^{2}+t^{3},\ldots$ = the number of hex tree by edges and left edges (\seqnum{A126182}).
\item $F(1, 1+t; ;z)\to  1,2+t, 5+ 4t+t^{2}, 14+ 15t+ 6t^{2} +t^{3},\ldots$ = the number of hex tree by edges and median children (\seqnum{A126181}).
\item $F(;2 ;z)\to 1,3,10,35,\ldots$ = the number of times a fixed Dyck word of length 2k occurs in all Dyck words of length $2n+2k$ = $\binom{2n+1}{n+1}$, (\seqnum{A001700}).
\item $F(;1+t;z)\to 1,2+t,5+4t+t^{2},14+14t+6t^{2}+t^{3},\ldots$ = the number of leaves at level k+1 in all ordered trees with n+1 edges (\seqnum{A039598}).
\item $F(0,1;1+t;z)\to 1,1+t,1+2t+t^{3},1+3t+3t^{2}+1,\ldots$ = Pascal's triangle (\seqnum{A007318}).
\item $F(0,1+t;1;z)\to 1,1+t,1+2t+t^{3},1+3t+3t^{2}+1,\ldots$ = Pascal's triangle (\seqnum{A007318}).
\item $F(1,0;t;z)\to 1,1,1+t,2+2t,4+4t+t^{2}, \ldots$ = the number of Motzkin paths by height of final descent (\seqnum{A098979}).
\end{itemize}
Note that in all these cases Corollary \ref{corollary:recursion} can be use to find an quadratic equation for the generating function.  Conversely, we have the following result. 
\begin{corollary}[Corollary 3.2, \cite{DFR2005}]
If a generating function $f$ satisfies an equation of the form
\[x_{0}y_{1}z^{2} f^{2}-\left(1-(x_{1}y_{1}+x_{0})z\right)f+1=0,\]
then $f=F(x_{0},x_{1};y_{1};z)$. 
\end{corollary}

%%%
%%%
%%%

\bibliographystyle{abbrv}

\end{document}